\documentclass[12pt]{article}

\usepackage{amsmath}
\usepackage{amsfonts}
\usepackage{amssymb}
\usepackage{graphicx, color}
\usepackage{hyperref}
\usepackage{epsfig}
\usepackage{epstopdf} 

\definecolor{darkgreen}{rgb}{0,0.55,0}

\newtheorem{proposition}{Proposition}[section]
\newtheorem{theorem}{Theorem}[section]
\newtheorem{lemma}[theorem]{Lemma}
\newtheorem{corollary}[theorem]{Corollary}
\newtheorem{remark}[theorem]{Remark}

\newtheorem{definition}{Definition}

\DeclareSymbolFont{AMSb}{U}{msb}{m}{n}
\DeclareMathSymbol{\N}{\mathbin}{AMSb}{"4E}
\DeclareMathSymbol{\Z}{\mathbin}{AMSb}{"5A}
\DeclareMathSymbol{\R}{\mathbin}{AMSb}{"52}
\DeclareMathSymbol{\Q}{\mathbin}{AMSb}{"51}
\DeclareMathSymbol{\I}{\mathbin}{AMSb}{"49}
\DeclareMathSymbol{\C}{\mathbin}{AMSb}{"43}

\begin{document}

\title{Determining both the source of a wave and its speed in a medium from boundary measurements }
\author{{Christina Knox\footnote{Department of Mathematics, University of California, Riverside, California, USA. E-mail: knox@math.ucr.edu.  }
\qquad Amir Moradifam \footnote{Department of Mathematics, University of California, Riverside, California, USA. E-mail: amirm@math.ucr.edu. }}}
\date{}

\smallbreak \maketitle

\begin{abstract}
We study the inverse problem of determining both the source of a wave and its speed inside a medium from measurements of the solution of the wave equation on the boundary.  This problem arises in photoacoustic and thermoacoustic tomography, and has important applications in medical imaging.  We prove that if the solutions of the wave equation with the source and sound speeds $(f_1,c_1)$ and $(f_2,c_2)$ agree on the boundary of a bounded region $\Omega$, then 

\[ \int_{\Omega}(c_2^{-2}-c_1^{-2})\varphi dy=0,\]
for every harmonic function $\varphi \in C(\bar{\Omega})$, which holds without any knowledge of the source.  We also show that if the wave speed $c$ is known and only assumed to be bounded then, under a natural admissibility assumption, the source of the wave can be uniquely determined from boundary measurements. 
\end{abstract}

\section{Introduction} 
Consider the wave equation 
\begin{eqnarray}\label{waveEq}
\left\{ \begin{array}{ll}
u_{tt}-c^2\Delta u=0 &\text{in } \R^3 \times \R_{+}\\
u(x,0)=f(x), \ \ u_t(x,0)=0 &\text{for } x\in \R^3,
\end{array} \right.
\end{eqnarray}
where $\Omega \subset \R^3$ is a simply connected bounded region and the wave speed $c \in L^{\infty}(\R^3)$ satisfies $c\geq c_0>0$ in $\Omega$, and supp$(1-c) \subset \subset \Omega$.  The function $f(x)\in L^{\infty}(\R^3)$ represents the source of the wave and is assumed to be compactly supported in $\Omega$, i.e. supp$(f) \subset \subset \Omega$.  
In this paper we study the inverse problem of recovering the pair $(f,c)$ from the measurements of the solution of the wave equation on $\partial \Omega$ given by the measurement operator
\begin{equation}
\Lambda_{f,c}(x,t)=u(x,t), \ \ \ \  (x,t) \in \partial \Omega \times \R_{+}. 
\end{equation}
This problem naturally arises in thermoacoustic (TAT) and photoacoustic (PAT) tomography, both of which have significant potential in clinical applications and biology \cite{HSBP, KKRK1, KRK, KRK,Wang, XW}. 

Unique determination of the source function $f$ and the wave speed $c$ has been studied by many authors and several interesting results have been obtained. However, most of the results in the literature have been concerned with determination of $f$ or $c$ from the knowledge of $\Lambda_{f,c}$ under the assumption that the other one is known. When the sound speed is known, smooth, and non-trapping, then the source $f$ can be uniquely recovered   \cite{AKK, FR, H, HKN, KK, SU2, SU1}. The recovery of a smooth sound speed when $f$ is known is studied in \cite{SU1}.  A stability estimate is obtained in \cite{OU} for the recovery of the source $f$ when there is a small error in the variable sound speed.  In practice the sound speed inside the medium is often unknown \cite{JW}. It has been observed that even replacing a sound speed with small variation by its average value can significantly distort the reconstruction of $f$ \cite{HKN}. One suggested solution is to additionally perform an ultrasonic transmission tomography (UTT) to recover the sound speed \cite{JW}. Thus from both a theoretical and practical point of view it would be advantageous to know whether  both the sound speed and the source term can be uniquely recovered from $\Lambda_{f,c}$. This is an open problem that we tackle in this paper. 

The first result for the recovery of a unknown sound speed was in \cite{HKN} where the authors proved that a constant sound speed can be uniquely recovered using range support conditions. It is shown via a connection to the transmission eigenvalue problem in \cite{FH} that if the sound speed is radial then both $f$ and $c$ can be recovered uniquely. In \cite{LU} H. Liu and G. Uhlmann showed that under additional assumptions on the wave speed and the source term both can be uniquely recovered simultaniously. In \cite{SU3} the authors proved that when both the sound speed and source are unknown the linearized problem is unstable.

In this paper, inspired by Liu and Uhlmann's approach in \cite{LU},  We prove that if  $c^{-2}$ is harmonic in $\omega \subset \R^3$ and identically 1 on $\omega^c$, where $\omega$ is a simply connected region, then a non-trapping wave speed $c$ can be uniquely determined from the solution of the wave equation on boundary of $\Omega \supset \supset \omega$ without the knowledge of the source.  We also show that if the wave speed $c$ is known and only assumed to be bounded then, under an admissibility assumption (see Definition 2 below), the source of the wave can be uniquely determined from boundary measurements. Indeed we prove Theorem \ref{mainTheorem} and  Theorem \ref{sourceUniqueness} below. To state the main results let us set the stage with two definitions.

 Let $u(x,t)$ be the solution of the wave equation \eqref{waveEq} and for $(x,k)\in \R^3 \times \R_{+}$ define the temporal Fourier transform of the function $u(x,t)$ by 
\begin{equation}\label{temporalFourier}
\hat{u}(x,k):=\frac{1}{2\pi} \int_{0}^{\infty}u(x,t)e^{ikt}dt. 
\end{equation}
Following \cite{LU}, let us define an admissible pair for a source and wave speed. 
\begin{definition} \label{admissiblality}
Let $0<c_0<c\in L^{\infty}(\Omega)$ and $f\in L^{\infty}(\Omega)$. We say that the pair $(f,c)$ is admissible if there exists $\epsilon>0$ such that $\hat{u}(x,k)\in H^1_{loc}(\R^3)$ for all $k\in (0,\epsilon)$. 
\end{definition}

Indeed if $u(x,t)$ decays fast enough in time such that $\hat{u}(x,k) \in H^1_{loc} (\R^3)$, then $(f,c)$ is admissible.  Note that the admissibility assumption above is a weak form of the non-trapping assumption on the sound speed. As pointed out in \cite{LU}, if the sound speed $c$ is smooth and non-trapping then $(f,c)$ is admissible, but not vice versa. See \cite{FH,SU1} and the references cited therein for more details. Throughout the paper we shall assume that the pair $(f,c)$ is admissible. 

Now we are ready to state the main results.

%For $(x,\xi)\in \R^{2n}_{x,\xi}$ and the Hamiltonian $H=\frac{1}{2}c^2(x)|\xi|^2$, consider the Hamiltonian system

%\begin{eqnarray}\label{Hamiltonian}
%\left\{ \begin{array}{ll}
%x'_t=\frac{\partial H}{\partial \xi}=c^2\xi \\
%\xi'_t=-\frac{\partial H}{\partial x}=-c(x)\nabla (c(x)) |\xi|^2,\\
%x|_{t=0}=x_0, \xi|_{t=0}=\xi_0.
%\end{array} \right.
%\end{eqnarray}
%The solution $(x(t),\xi(t))$ is called a bicharacteristic and $x(t)$ is called a ray. 
%\begin{definition}\label{nontrapping}
%We say that the sound speed $c$ is non-trapping if all rays with $\xi_0\neq 0$ tend to infinity as $t\rightarrow \infty$. 

%\end{comment}

\begin{theorem}\label{mainTheorem}
Let $(f_1,c_1)$ and $(f_2,c_2)$ be two admissible pairs such that 
\[\Lambda_{f_1,c_1}(x,t)=\Lambda_{f_2,c_2}(x,t) \ \ \ \ \forall (x,t)\in \partial \Omega \times \R_{+}.\]
Then  
\[C^*:=\int_{\Omega}c_1^{-2}(x)f_1(x)dx=\int_{\Omega}c_2^{-2}(x)f_2(x)dx.\]
If $C^*\neq 0$, then

\begin{equation}\label{mainCond}
\int_{\Omega}(c_2^{-2}-c_1^{-2})\varphi dy=0, \ \ \hbox{for all harmonic functions }  \varphi.
\end{equation}
In particular, if $(c_2^{-2}-c_1^{-2})$ is harmonic in a simply connected region $\omega \subset \subset \Omega$ and identically zero on $\omega^c$, then $c_1 \equiv c_2$ in $\Omega$. 
\end{theorem}

Note that if \eqref{mainCond} holds for every smooth function $\varphi$, then it will imply $c_1=c_2$. Indeed Theorem \eqref{mainTheorem} is an step towards simultaneous recovery of the sound speed and the source in the wave equation. It should be compared to the stronger results in \cite{SU1} where the authors prove uniqueness of the sound speed from the knowledge of the source $f$ and under the assumption that the domain $\Omega$ is foliated by strictly convex hypersurfaces with respect to the Riemannian metric $g=c^{-2}dx$.

\begin{theorem}\label{sourceUniqueness} 
Let $0<c_0<c \in L^{\infty}(\Omega)$ for some $c_0>0$, and $(f_1,c)$ and $(f_2,c)$ be admissible pairs. If  
\[\Lambda_{f_1,c}(x,t)=\Lambda_{f_2,c}(x,t) \ \ \ \ \forall (x,t)\in \partial \Omega \times \R_{+},\]
then  $f_1=f_2$. 
\end{theorem}

Theorem \ref{sourceUniqueness} should be compared to the uniqueness results in \cite{AKK} and \cite{SU2} where the authors assume that the sound speed is smooth. Models with discontinuous sound speed arise in thermoacoustic and photoacoustic tomography in order to understand the effect of sudden change of the sound speed in the skull in imaging of the human brain \cite{SU4}. The results in \cite{SU4} assume that the sound speed is smooth but allow for jumps across smooth surfaces. In \cite{SU2} and \cite{SU4} the authors only require the knowledge $u(x,t)$ on $\partial \Omega \times [0,T]$ for some finite time $T>0$ while Theorem \ref{sourceUniqueness} assumes the knwoeledge of the solution on $\partial \Omega \times (0,\infty)$. Explicit reconstruction formulas for such problems have been also developed in \cite{AKK, SU2, SU4}.

\vspace{-.3cm}
\section{Uniqueness of the wave speed}
In this section we present the proof of Theorem \ref{mainTheorem}. Let us first develop a few basic facts about solutions of the wave equation \eqref{waveEq} and gather some known results which will be used in our proofs. 

The temporal Fourier transform $\hat{u}$, define in \eqref{temporalFourier}, satisfies the elliptic partial differential equation
\begin{equation}\label{u-hatPDE}
\Delta \hat{u}(x,k)+\frac{k^2}{c^2(x)}\hat{u}(x,k)=\frac{ik}{2\pi} \frac{f(x)}{c^2(x)}, \ \ (x,k)\in \R^3 \times (0,\epsilon),
\end{equation}
which is well-posed under the classical Sommerfeld radiation condition
\begin{equation}\label{sommerfeld}
\lim_{|x|\rightarrow \infty} |x| \left( \frac{\partial \hat{u}(x,k)}{\partial |x|}-ik\hat{u}(x,k)\right)=0. 
\end{equation}
Note also that under the temporal Fourier transform the measurement operator becomes
\begin{equation}
\hat{\Lambda}_{f,c}(x,k)=\hat{u}(x,k), \ \ \ \ (x,k) \in \R^3 \times (0,\epsilon). \\ \\
\end{equation}
We shall need the following lemma proved by Liu and Uhlamm in \cite{LU}. 
\begin{lemma} [\cite{LU}] \label{LiuUhlmannLemma} Let $\hat{u}(x,k)\in H^1_{loc}(\R^3)$ be the solution to \eqref{u-hatPDE}-\eqref{sommerfeld}. Then $\hat{u}(x,k)$ is uniquely given by the following integral equation 
\begin{equation}\label{RepresentationFormula}
\hat{u}(x,k)=k^2\int_{\R^3}\left( c^{-2}(y)-1\right)\hat{u}(y,k) \Phi(x-y)dy-\frac{ik}{2\pi}\int_{\R^3}\frac{f(y)}{c^{2}(y)}\Phi(x-y)dy,  \ \ x\in \R^3.
\end{equation}
Moreover, as $k\rightarrow 0$, we have 
\begin{eqnarray}\label{firstTwoTerms}
\hat{u}(x,k)=-\frac{ik}{2\pi}\int_{\Omega}\frac{f(y)}{c^{2}(y)} \Phi_0(x-y)dy+\frac{k^2}{8\pi^2}\int_{\Omega}\frac{f(y)}{c^{2}(y)}dy+ \mathcal{O}(k^3).
\end{eqnarray}
 
\end{lemma}
Here 
\[\Phi(x):=\frac{e^{ik|x|}}{4\pi|x|} \ \ \hbox{for}\ \ |x|\neq 0\]
is the fundamental solution of $-\Delta-k^2$ and $\Phi_0$ is the fundamental solution of $-\Delta$. 

Define the space
\[\mathcal{A}:=\{v\in L^2(\Omega): \int_{\Omega}g \varphi dx=0 \ \ \hbox{for all harmonic functions }  \varphi\in C(\bar{\Omega)} \}.\]
We shall frequently use the following two lemmas. 

\begin{lemma}\label{PoissonLemma}
Let $g\in L^2(\Omega)$ and suppose $w \in H^1(\Omega)$ satisfies 
\begin{eqnarray}\label{Poisson}
\left\{ \begin{array}{ll}
-\Delta w=g &\text{in } \Omega\\
w=0 \ \ &\text{on } \partial \Omega.
\end{array} \right.
\end{eqnarray}
Then $\frac{\partial w}{\partial \nu}=0$ on $\partial \Omega$ if and only if $g\in \mathcal{A}$. 
\end{lemma}
{\bf Proof.} Let $\varphi \in C(\bar{\Omega})$ be harmonic in $\Omega$ and $w \in H^2(\Omega)$ be the solution of \eqref{Poisson}. By integration by parts one can show that
\[-\int_{\partial\Omega}\varphi(y) \frac{\partial w}{\partial \nu}dS=\int_{\Omega}g(y)\varphi(y)dy.\]
Hence $\frac{\partial w}{\partial \nu}=0$ on $\partial \Omega$ if and only if $g\in \mathcal{A}$. \hfill $\Box$

\begin{lemma}\label{NewtonianPotLemma}
For $g \in L^{\infty}(\R^3)$ with compact support in $\Omega$ define 
\begin{equation}\label{DeltaInverse}
w(x):=\int_{\R^3}g(y)\Phi_0(x-y)dy=\int_{\Omega}g(y)\Phi_0(x-y)dy.
\end{equation}
Then 
\[-\Delta w=g \ \ \hbox{in} \ \ \R^3\]
in the weak sense. Moreover if $w=0$ on $\Omega^c$, then for any harmonic function $\varphi$ on $\R^3$ we have
\[\int_{\R^3}g(y)\varphi(y) dy=\int_{\Omega}g(y)\varphi(y)dy=0.\]
\end{lemma}
{\bf Proof.}  Since  $w=0$ and $\frac{\partial w}{\partial \nu}=0$, the result follows from Lemma \ref{PoissonLemma}. \hfill $\square$ 

For $g \in L^2(\Omega)$ we will denote the solution of \eqref{Poisson} by $\Delta^{-1}(g)$. Note that if  $g \in L^{\infty}(\R^3)$ has compact support in $\Omega$  and $w$ defined by \eqref{DeltaInverse} vanishes on $\Omega^c$, then $w=\Delta^{-1}(g)$.  \\ \\ 
{\bf Proof of Theorem \ref{mainTheorem}.} By Lemma \ref{LiuUhlmannLemma} we have 
\begin{eqnarray}\label{firstTwoTerms}
\hat{u}(x,t)=-\frac{ik}{2\pi}\int_{\Omega}\frac{f(y)}{c^{2}(y)}dy \Phi_0(x-y)dy+\frac{k^2}{8\pi^2}\int_{\Omega}\frac{f(y)}{c^{2}(y)}dy+ \mathcal{O}(k^3). 
\end{eqnarray}
For $i=1,2$ define 
\[w_i(x):=\lim _{k\rightarrow 0} \frac{\hat{u}_i(x,k)}{k}=-\frac{i}{2\pi} \int_{\Omega}\frac{f_i(y)}{c_i^{2}(y)} \Phi_0(x-y)=\frac{i}{2\pi}\Delta^{-1}(\frac{f_i(y)}{c_i^{2}(y)}).\]
Then $w:=w_2-w_1$ satisfies
\[\Delta w= \frac{i}{2\pi} \left(\frac{f_2(y)}{c_2^{2}(y)}-\frac{f_1(y)}{c_1^{2}(y)}  \right),\]
and $w=\frac{\partial w}{\partial \nu}=0$ on $\partial \Omega$. Therefore it follows from Lemma \ref{PoissonLemma} that 
\begin{equation}\label{orthogonalHarmonic0}
\int_{\Omega}\left( \frac{f_2(y)}{c_2^{2}(y)}-\frac{f_1(y)}{c_1^{2}(y)} \right) \varphi(y)dy =0,
\end{equation}
for every harmonic function $\varphi$. On the other hand we have 
\begin{equation}\label{difference}
\Delta ( \hat{u}_2(x,k)-  \hat{u}_1(x,k))+\frac{k^2}{c_2^2(x)}\hat{u}_2(x,k)-\frac{k^2}{c_1^2(x)}\hat{u}_1(x,k)=-\frac{i k}{2\pi}\left( \frac{f_2(x)}{c_2^{2}(x)}-\frac{f_1(x)}{c_1^{2}(x)} \right),
\end{equation}
for $ (x,k)\in \R^3 \times (0,\epsilon).$
Multiplying both sides of the above equation by a harmonic function $\varphi$, using \eqref{orthogonalHarmonic0} and the fact that $\hat{u}_2-\hat{u}_1\equiv 0$ on $\Omega^c$, and integrating by parts we get 
\begin{equation}\label{multiplyvarphi}
\frac{1}{k} \int_{\Omega} \left(  \frac{\hat{u}_2(y,k)}{c_2^2(y)}-\frac{\hat{u}_1(y,k)}{c_1^2(y)} \right) \varphi dy=0, \ \ \ \ \forall k\in (0,\epsilon). 
\end{equation}
Since $\hat{u}_2(x,k)=\hat{u}_1(x,k)$ for all $(x,k)\in \partial \Omega \times (0,\epsilon)$, it follows from \eqref{firstTwoTerms} that 
\begin{equation}\label{independent}
\int_{\Omega}c_1^{-2}f_1dy=\int_{\Omega} c_2^{-2}f_2dy.
\end{equation}
Combining this with \eqref{multiplyvarphi} and Lemma \ref{LiuUhlmannLemma} we have
\begin{eqnarray}
&& \frac{i}{2\pi}\int_{\Omega} \left( \frac{(\int_{\Omega}c_2^{-2}(z)f_2(z)\Phi_0(y-z)dz}{c_2^2(y)}- \frac{\int_{\Omega}c_1^{-2}(z)f_1(z)\Phi_0(y-z)dz}{c_2^2(y)} \right)\varphi (y)dy \nonumber \\ && -\frac{k}{8\pi^2}\int_{\Omega} c_1^{-2}(y)f_1(y) dy \int_{\Omega}(c_2^{-2}(y)-c_1^{-2}(y))\varphi (y)dy+\mathcal{O}(k^2)=0,
\end{eqnarray}
for all $k\in (0,\epsilon)$. Thus 
\begin{equation}
\int_{\Omega} \left( \frac{\int_{\Omega}c_2^{-2}(z)f_2(z)\Phi_0(y-z)dz}{c_2^2(y)}- \frac{\int_{\Omega}c_1^{-2}(z)f_1(z)\Phi_0(y-z)dz}{c_2^2(y)} \right)\varphi (y)dy
\end{equation}
and
\begin{equation}\label{c_1-c_2}
\int_{\Omega}(c_2^{-2}-c_1^{-2})\varphi dy=0
\end{equation}
for any harmonic function $\varphi$, provided $C^*\neq 0$. Finally note that if $(c_2^{-2}-c_1^{-2})$ is harmonic, then letting $\varphi=(c_2^{-2}-c_1^{-2})$ in \eqref{c_1-c_2} implies $c_1 \equiv c_2$.  \hfill $\Box$

\begin{corollary}\label{cor} If $c_1 \geq c_2$ or $c_1 \leq c_2$ in $\Omega$, then $c_1\equiv c_2$. 
\end{corollary}
{\bf Proof.} Let $\varphi\equiv 1$ in $\Omega$. Then by Theorem \ref{mainTheorem} we have 
\[\int_{\Omega} c_2^{-2}(x)dx= \int_{\Omega}c_1^{-2}(x)dx,\]
and the result follows immediately. \hfill $\Box$\\ \\
See \cite{FH} for another proof of Corollary \ref{cor}. 

\begin{remark}
Note that $\hat{u}_1$ and $\hat{u}_2$ can be represented as 
\[\hat{u}_i(y,k)=\sum \limits_{n=1}^{\infty} p_n^i(y)k^n, \ \ \hbox{for}\ \ (y,k)\in \R^3 \times (0,\epsilon), \ \  i=1,2\]
(see Proposition \ref{formulaforpn}).
If $\hat{u}_1 (y,k)= \hat{u}_2(y,k)$ for $(y,k)\in \partial \Omega \times (0,\epsilon)$, then 
\begin{equation}\label{a_n}
p_n^1(y)=p_n^2(y), \ \ \forall y \in \partial \Omega \ \ \hbox{and}\ \  \forall n \in \N. 
\end{equation}
The proof of Theorem \ref{mainTheorem} only uses \eqref{a_n} for $n=1, 2$. One may expect to combine the condition \eqref{a_n} for $n\geq 2$ and the equality \eqref{multiplyvarphi} to prove that the sound speed can be uniquely determined from boundary measurements,  without any major assumptions on the sound speed. This strongly suggests that $\Lambda_{f,c}$ could uniquely determine both the source and the sound speed in general. The authors believe that the measurement operator $\Lambda_{f,c}$ can uniquely determine both the sound speed $c$ and the source $f$. However, the higher order terms turn out to be complicated and the authors were not able to prove the complete result. 

\end{remark}

\section{Uniqueness of the source}
In this section we prove that if $c_1=c_2=c\in L^{\infty}$, then the source function $f$ can be uniquely recovered from the knowledge of $\Lambda_{f,c}(x,t)$ on $ \partial \Omega \times \R_+$. Throughout this section we shall assume that $c_1=c_2=c \in L^{\infty}(\R^3)$.\\

By Lemma \ref{LiuUhlmannLemma}, $\hat{u_2}-\hat{u_1}$ satisfies the following integral equation 
\begin{eqnarray}\label{inteqndiff}
(\hat{u_2}-\hat{u_1})(x,k)&=&k^2\int_{\Omega}(c^{-2}(y)-1)(\hat{u_2}- \hat{u_1})(y,k)\Phi(x-y)dy  \nonumber \\ 
&& -\frac{ik}{2\pi}\int_{\Omega}c^{-2}(y)(f_2-f_1)(y)\Phi(x-y)dy.
\end{eqnarray}
We shall need the following lemma. 
\begin{lemma}\label{IbyP} Let $g$ have compact support in $\Omega$ and $g \in L^{\infty}(\R^3)$. Suppose that $w$ defined by \eqref{DeltaInverse} vanishes on $\Omega^c$. Then 
\begin{equation}
\int_{\Omega}g(y)|x-y|^ndy=-n(n+1)\int_{\Omega}\Delta^{-1}(g)(y)|x-y|^{n-2}dy
\end{equation}
for all $x \in \R^3$.
\end{lemma}
{\bf Proof.} The proof follows from Lemma \ref{NewtonianPotLemma} and integration by parts. \hfill $\Box$\\ 

For every $g \in L^{\infty}(\R^3)$ with compact support in $\Omega$ define 
\[Lg(x):=\Delta^{-1}(c^{-2}g)(x), \ \ x \in \R^3.\]

\begin{proposition}\label{formulaforpn} For every $n\in \N$ there exists functions $p_m(x)$, $m=1,2, . . . , n$, such that 
\begin{equation}\label{firstStatement}
(\hat{u_2}-\hat{u_1})(x,k)=\sum \limits_{m=1}^{n}p_m(x)k^{j}+\mathcal{O}(k^{n+1}),
\end{equation}
as $k\rightarrow 0$. Moreover if $u_2(x,k)-u_1(x,k)=0$ for all $x \in \Omega^c$,
then
\begin{equation}\label{p_n}
p_n(x)=
\begin{cases}
      \frac{-i}{2\pi}L^{\frac{n+1}{2}}(f_2-f_1)(x) & \textnormal{ if } n  \textnormal{ is odd} \\
     0 & \textnormal{ if } n  \textnormal{ is even.}
    \end{cases}
\end{equation}
\end{proposition}
{\bf Proof.} By Lemma \ref{LiuUhlmannLemma}, \eqref{firstStatement} holds for $n=1,2$. Suppose it holds for all $j \leq n$. Then there exists functions $p_m(x)$, $m=1,2, . . . , n$,  such that
\[(\hat{u_2}-\hat{u_1})(x)=\sum \limits_{m=1}^{n}p_m(x)k^{j}+\mathcal{O}(k^{n+1}) \ \ \hbox{as} \ \ k \rightarrow 0.\]
Plugging this expression for $\hat{u_2}-\hat{u_1}$ into equation \eqref{inteqndiff} and expanding $\Phi$ we find that
\[(\hat{u_2}-\hat{u_1})(x)=\sum \limits_{m=1}^{n}p_m(x)k^{j}+p_{n+2}(x)k^{n+2}+\mathcal{O}(k^{n+3}) \ \ \hbox{as} \ \ k \rightarrow 0,\]
where
\begin{align*}
p_{n+2}(x) &= \sum_{m=0}^{n}\frac{i^m}{4\pi m!}\int_{\Omega} (c^{-2}-1)(y)p_{n-m}(y)|x-y|^{m-1}dy \\
& \ \ \ \ - \frac{i^{n+2}}{8\pi^2(n+1)}\int_{\Omega} c^{-2}(y)(f_2-f_1)(y)|x-y|^n dy.
\end{align*}
To prove \eqref{p_n} we proceed by strong induction. First notice that $p_m\equiv 0$ on $\Omega^c$ for all $m\in \N$. By Lemma \ref{LiuUhlmannLemma}, \eqref{p_n} holds when $n=0,1,2$. Suppose \eqref{p_n} holds for all $j \leq n+1$. First assume that $n$ is odd. Using the integral equation \eqref{inteqndiff} and the induction hypothesis we compute that
\begin{align*}
p_{n+2}(x) &= \sum_{\substack{m=0 \\ \textnormal{m even}}}^{n-1}\frac{i^{m+1}}{8\pi^2 m!}\int_{\Omega} (1-c^{-2})(y)L^{\frac{n-m+1}{2}}(f_2-f_1)(y)|x-y|^{m-1}dy\\ 
& \ \ \ \ -\frac{i^{n+2}}{8\pi^2(n+1)}\int_{\Omega} c^{-2}(y)(f_2-f_1)(y)|x-y|^n dy.
\end{align*}
For even $m$  with $m \leq n-1$ define 
\[q_m(x):=\frac{i^{m+1}}{8\pi^2 m!}\int_{\Omega} (1-c^{-2})(y)L^{\frac{n-m+1}{2}}(f_2-f_1)(y)|x-y|^{m-1}dy,\]
and 
\[r(x):=-\frac{i^{n+2}}{8\pi^2(n+1)!}\int_{\Omega} c^{-2}(y)(f_2-f_1)(y)|x-y|^n dy.\]
Then 
\[p_{n+2}(x) = \sum_{\substack{m=0 \\ \textnormal{m even}}}^{n-1}q_m(x)+r(x).\]
It follows from the induction hypothesis and Lemma \ref{IbyP} that
\begin{align*}
q_2(x)&=\frac{-i}{8\pi^2 2!}\int_{\Omega} (1-c^{-2})(y)L^{\frac{n-1}{2}}(f_2-f_1)(y)|x-y|dy\\
&=\frac{-i}{8\pi^2 2!}\int_{\Omega}L^{\frac{n-1}{2}}(f_2-f_1)(y)|x-y|dy\\
& \ \ \ \ +\frac{i}{8\pi^2 2!}\int_{\Omega}c^{-2}(y)L^{\frac{n-1}{2}}(f_2-f_1)(y)|x-y|dy\\
&=\frac{-i}{8\pi^2 2!}\int_{\Omega}L^{\frac{n-1}{2}}(f_2-f_1)(y)|x-y|dy\\
& \ \ \ \ -\frac{i}{8\pi^2}\int_{\Omega}\Delta^{-1}(c^{-2}(y)L^{\frac{n-1}{2}}(f_2-f_1))(y)|x-y|^{-1}dy\\
&=\frac{-i}{8\pi^2 2!}\int_{\Omega}L^{\frac{n-1}{2}}(f_2-f_1)(y)|x-y|dy-\frac{i}{8\pi^2}\int_{\Omega}L^{\frac{n+1}{2}}(f_2-f_1)(y)|x-y|^{-1}dy.
\end{align*}
Thus we have 
\begin{align*}
(q_0+q_2)(x)&=\frac{-i}{8\pi^2}\int_{\Omega}c^{-2}(y)L^{\frac{n+1}{2}}(f_2-f_1)(y)|x-y|^{-1}dy \\ 
& \ \ \ \ -\frac{i}{8\pi^2 2!}\int_{\Omega}L^{\frac{n-1}{2}}(f_2-f_1)(y)|x-y|dy\\
&=\frac{-i}{2\pi}\int_{\Omega}c^{-2}(y)L^{\frac{n+1}{2}}(f_2-f_1)(y)\Phi_0(x-y)dy\\
&-\frac{i}{8\pi^22!}\int_{\Omega}L^{\frac{n-1}{2}}(f_2-f_1)(y)|x-y|dy\\
&=-\frac{i}{2\pi}L^{\frac{n+3}{2}}(f_2-f_1)(x)-\frac{i}{8\pi^22!}\int_{\Omega}L^{\frac{n-1}{2}}(f_2-f_1)(y)|x-y|dy.
\end{align*}
Similarly by Lemma \ref{IbyP} we get 
\begin{eqnarray*}
q_4(x) &=& \frac{i}{8\pi^2 4!}\int_{\Omega} (1-c^{-2})(y)L^{\frac{n-3}{2}}(f_2-f_1)(y)|x-y|^{3}dy\\
&=& \frac{i}{8\pi^2 4!}\int_{\Omega}L^{\frac{n-3}{2}}(f_2-f_1)(y)|x-y|^3dy+\frac{i}{8\pi^2 2!}\int_{\Omega}L^{\frac{n-1}{2}}(f_2-f_1)(y)|x-y|dy.
\end{eqnarray*}
Hence
\[(q_0+q_2+q_4)(x)=-\frac{i}{2\pi}L^{\frac{n+3}{2}}(f_2-f_1)(x)+\frac{i}{8\pi^2 4!}\int_{\Omega}L^{\frac{n-3}{2}}(f_2-f_1)(y)|x-y|^3dy.\]
We can continue this process in general. Let $m$ be even with $m \leq n-3$ and suppose
\begin{align*}
(q_0+q_2+...+q_{m-2}+q_m)(x)& = -\frac{i}{2\pi}L^{\frac{n+3}{2}}(f_2 - f_1)(x)+\frac{i^{m+1}}{8\pi^2 m!}L^{\frac{n-m+1}{2}}(f_2-f_1)(y)|x-y|^{m-1}dy.
\end{align*}
Then by Lemma \ref{IbyP}
\begin{eqnarray*}
q_{m+2}(x) &=& \frac{i^{m+3}}{8\pi^2 (m+2)!}\int_{\Omega} (1-c^{-2})(y)L^{\frac{n-m-1}{2}}(f_2-f_1)(y)|x-y|^{m+1}dy\\
&=& \frac{i^{m+3}}{8\pi^2 (m+2)!}\int_{\Omega}L^{\frac{n-m-1}{2}}(f_2-f_1)(y)|x-y|^{m+1}dy\\
&& +\frac{i^{m+3}}{8\pi^2 m!}\int_{\Omega}L^{\frac{n-m+1}{2}}(f_2-f_1)(y)|x-y|^{m-1}dy.
\end{eqnarray*}
Noting that $i^{m+3}=-i^{m+1}$ we get that
\begin{eqnarray*}
(q_0+q_2+...+q_m+q_{m+2})(x)&=&-\frac{i}{2\pi}L^{\frac{n+3}{2}}(f_2-f_1)(x)\\
&&+\frac{i^{m+3}}{8\pi^2 (m+2)!}\int_{\Omega}L^{\frac{n-m-1}{2}}(f_2-f_1)(y)|x-y|^{m+1}dy.
\end{eqnarray*}
Repeating the above process until $m=n-1$ we obtain 
\begin{eqnarray*}
(q_0+q_2+...+q_{n-3}+q_{n-1})(x)&=&-\frac{i}{2\pi}L^{\frac{n+3}{2}}(f_2-f_1)(x)\\
&& +\frac{i^n}{8\pi^2 (n-1)!}\int_{\Omega}L(f_2-f_1)(y)|x-y|^{n-2}dy.
\end{eqnarray*}
In addition by Lemma \ref{IbyP}
\[r(x)=\frac{i^{n+2}}{8\pi^2(n-1)!}\int_{\Omega} L(f_2-f_1)(y)|x-y|^{n-2} dy.\]
Hence
\[p_{n+2}(x)=-\frac{i}{2\pi}L^{\frac{n+3}{2}}(f_2-f_1)(x).\]
This finishes the proof for the case that $n$ is odd. 

Now suppose $n$ is even.
With an argument similar to the one in the case when $n$ was odd we can show
\[p_{n+2}(x)=\frac{1}{8\pi^2}\int_{\Omega}c^{-2}(y)L^{\frac{n}{2}}(f_2-f_1)(y)dy.\]
By the induction hypothesis for $j=n+1$ we have 
\[p_{n+1}(x)=L^{\frac{n+2}{2}}(f_2-f_1)(x)=0, \ \ \forall x \in \Omega^c.\]
Hence 
\[\Delta^{-1}(c^{-2}L^{\frac{n}{2}}(f_2-f_1))=L^{\frac{n+2}{2}}(f_2-f_1)(x)=0 \ \ \forall x \in \Omega^c,\]
and it follows from Lemma \ref{NewtonianPotLemma} that 
\[\int_{\Omega}c^{-2}(y)L^{\frac{n}{2}}(f_2-f_1)(y) \varphi dy=0,\]
for every harmonic function $\varphi.$ Letting $\phi \equiv 1$ we get 
\[p_{n+2}(x)=\int_{\Omega}c^{-2}(y)L^{\frac{n}{2}}(f_2-f_1)(y)dy=0.\]
 Hence $p_{n+2}(x)=0$. \hfill $\square$

\begin{theorem}\label{PatternTheo}
Suppose $c_1=c_2=c$ and that $u_1,u_2$ be solutions of the wave equation \eqref{waveEq} with $u_1(x,0)=f_1$ and $u_2(x,0)=f_2$. If 
\[u_2(x,t)=u_1(x,t) \textnormal{ for all } (x,t) \in \Omega^c \times \R_+,\]
then 
\begin{equation}
\int_{\Omega} F_n\varphi c^{-2}dx=0, 
\end{equation}
for all $n\geq 0$ and all harmonic functions $\varphi$, where $F_0=f_2-f_1$ and 
\begin{equation}\label{F_n}
F_n=\Delta^{-1}(c^{-2} F_{n-1}), \ \ n\geq 1.
\end{equation}
\end{theorem}
{\bf Proof.} The proof follows directly from Proposition \ref{formulaforpn}, Lemma \ref{PoissonLemma}, and the observation that $F_n=p_n$. \hfill $\square$

\begin{lemma}\label{LemmaNoClusterCriticalPoints}
Let $g \in L^2(\Omega)$ satisfy 
\begin{equation}\label{lemma.orthog1}
\int_{\Omega} c^{-2} g \varphi dx=0,
\end{equation}
for all $\varphi \in C(\bar{\Omega})$ harmonic in $\Omega$. Suppose $c_0< c \in L^{\infty}(\Omega)$ for some $c_0>0$. If 
\begin{equation}\label{lemma.pde1}
 \lambda g=  \Delta^{-1} (c^{-2} g) \ \ \hbox{in} \ \ \Omega
\end{equation}
for some $\lambda>0$, then $g\equiv 0$. 
\end{lemma}
{\bf Proof.} It follows from Lemma \ref{PoissonLemma} that $g|_{\Omega}=\frac{\partial g}{\partial \nu}=0$. Note that
\[-\Delta g=\frac{c^{-2}}{\lambda}g \ \ \hbox{in} \ \ \Omega.\]
Since $c\in L^{\infty}(\Omega)$, by elliptic regularity $g \in C^{1,\alpha}(\Omega)$ for any $\alpha \in (0,1)$. Hence we can extend $g$ to a function $\tilde{g} \in C^{0,\alpha}(\R^3)$ by defining $\tilde{g}=0$ on $\Omega^c$. 
Let 
\[w(x):=\frac{1}{\lambda}\int_{\R^3}c^{-2}(y)\tilde{g}(y)\Phi_0(x-y)dy=\frac{1}{\lambda}\int_{\Omega}c^{-2}(y)g(y)\Phi_0(x-y)dy.\]
Since $c^{-2}\tilde{g} \in L^{\infty}(\R^3)$, it follows from elliptic regularity that $w \in C^{1,\alpha}(\R^3)$ and it satisfies
\[-\Delta w=\frac{c^{-2}}{\lambda}\tilde{g}.\]
Furthermore since $g \in \mathcal{A}$, $w=0$ on $\Omega^c$. Thus $w=\tilde{g}$ and hence $\tilde{g} \in C^{1,\alpha}(\R^3)$ solves
\[-\Delta \tilde{g}=\frac{c^{-2}}{\lambda}\tilde{g}, \ \ \hbox{in} \ \ \R^3\]
and $\tilde{g}=0$ on $\Omega^c$. Hence it follows from the unique continuation results in \cite{JK-UniqueContinuation} (see Theorem 6.3 ) that $\tilde{g}\equiv g \equiv 0$ in $\Omega$. 

  \hfill $\Box$

\begin{lemma}\label{Hilbertspace}
Let $\mathcal{H}$ be defined by 
\[\mathcal{H}:=\{v\in L^2(\Omega, c^{-2}dx):  c^{-2}v \in \mathcal{A} \ \ \hbox{and} \ \ c^{-2}L^n(v)\in \mathcal{A} \ \ \forall n\in \N\}.\]
Then $\mathcal{H}$ is a Hilbert space equipped with the inner product
\[\langle v,w\rangle|_{\mathcal{H}}=\int_{\Omega}vw c^{-2}dx.\]
\end{lemma}
\textbf{Proof.} The linearity of $L$ gives that $\mathcal{H}$ is a subspace. It remains to verify that $\mathcal{H}$ is closed. Suppose $v_n$ converges to $v$ in  $L^2(\Omega, c^{-2}dx)$ and $v_n \in \mathcal{H}$. It is easy to show that $c^{-2}v \in \mathcal{A}$. We claim that $w_n:=Lv_n$ converges to $w:=Lv$ in $L^2(\Omega, c^{-2}dx)$. Since 
\[-\Delta (w_n-w)=c^{-2}(v_n-v) \ \ \hbox{in} \ \ \Omega, \ \ w_n=0 \ \ \hbox{on} \ \ \partial \Omega,\]
we have
\[\parallel w_n-w \parallel_{H^2(\Omega)} \leq C \parallel c^{-2}(v_n-v) \parallel_{L^2(\Omega)} \rightarrow 0,\]
(See the Remark after Theorem 4 in section 6.3 of \cite{Evans}) and hence $w_n$ converges to $w$ in $L^2(\Omega, c^{-2}dx)$. Moreover, $c^{-2}L^n(w)=c^{-2}L^{n+1}(v)\in \mathcal{A}$ for all $n\geq 0$. Thus $w\in \mathcal{H}$ and the proof is complete.   \hfill  $\Box$
\begin{lemma}\label{PropsofL}
The linear operator $L: \mathcal{H} \rightarrow \mathcal{H}$ is a non-negative, self-adjoint, and compact operator.
\end{lemma}
{\bf Proof.} It follows from integration by parts that
\[\langle Lv,v\rangle_{\mathcal{H}}=\int_{\Omega} \Delta^{-1}(c^{-2}v) c^{-2}v=\int_{\Omega}|\nabla \Delta^{-1}(c^{-2}v)|^2 \geq 0,\]
and hence $L$ is non-negative. Similarly,
\[\langle Lv,w \rangle_{\mathcal{H}} = \int_{\Omega} \Delta^{-1}(c^{-2}v) wc^{-2}dx=\int_{\Omega}c^{-2}v \Delta^{-1}(c^{-2}w) dx=\langle v, Lw \rangle_{\mathcal{H}} ,\]
and hence $L$ is self-adjoint.
To show that $L$ is compact we need to to prove that $L(B_{\mathcal{H}})$ has compact closure in the strong topology, where $B_{\mathcal{H}}$ is the unit ball in $\mathcal{H}$ (see \cite{Brezis}). Let $v_n \in B_{\mathcal{H}}$. We need to show that $\{w_n\}:=\{ L(v_n) \}$ has a subsequence that converges in $L^2(\Omega,c^{-2}dx)$. Since
\[-\Delta w_n=c^{-2}v_n, \ \ w_n=0 \ \ \hbox{on} \ \ \partial \Omega,\]
we have
\[\parallel w_n \parallel_{H^2(\Omega)} \leq C \parallel c^{-2}v_n \parallel_{L^2(\Omega)} \leq C.\]
Thus $w_n$ is bounded in $H^2(\Omega)$ and hence $w_n$ has a subsequence, denoted by $w_n$ again, that converges weakly in $H^2(\Omega)$. Therefore $w_n$ converges strongly in $L^2(\Omega, c^{-2}dx)$ to some $w \in L^2(\Omega, c^{-2}dx)$ and thus $L$ is compact. \hfill $\square$
\begin{proposition}\label{SpectrumProp}
Let $F_n$ be defined by \eqref{F_n} and suppose 
\begin{equation}\label{PatternOdd}
\int_{\Omega} F_n\varphi c^{-2}dx=0, 
\end{equation}
for all $n\geq 0$ and all harmonic functions $\varphi \in C(\Omega)$. Then $F_0\equiv 0$ in $\Omega$.\\ 
\end{proposition}
{\bf Proof.} It follows from Lemma \ref{Hilbertspace} and  Lemma \ref{PropsofL} that $L$ has an orthonormal basis of eigenfunctions $e_n\in \mathcal{H}$ with corresponding eigenvalues $\lambda_n \geq 0$, where $\lambda_n \geq \lambda_{n+1}$ and $\lambda_n \rightarrow 0$ as $n\rightarrow \infty$.  Suppose $F_0\not \equiv 0$. Then there exists constants $\alpha_j \in \R$ such that 
\begin{equation}
F_0=\sum \limits_{j=1}^{\infty} \alpha_j e_j. 
\end{equation}
Let $\varphi \in C(\overline{\Omega})$ be harmonic in $\Omega$. Then \eqref{PatternOdd} implies 
\begin{eqnarray}\label{equality2}
\int_{\Omega} \left( \sum \limits_{j=1}^{\infty} \lambda_j^{n}\alpha_j e_j(x) \right)\varphi c^{-2} dx=0, \ \ \forall n\geq 0. 
\end{eqnarray}
Now let 
\[\lambda_{*}=\max \{\lambda_j:  \ \ \alpha_j\neq 0\}.\]
Dividing equality \eqref{equality2} by $\lambda_*^{n}$ yields
\[\int_{\Omega} \left( \sum \limits_{j=1}^m \alpha_j e_j \right)\varphi c^{-2} dx+\int_{\Omega} \left( \sum \limits_{j=m+1}^{\infty} \left(\frac{\lambda_j}{\lambda_{*}}\right)^{n}\alpha_j e_j \right)\varphi c^{-2} dx=0, \ \ \forall n\geq 0,\]
where $L(e_j)=\lambda_* e_j$, $j=1,2,... m$. 
We observe that 
\begin{eqnarray*}
\left \vert \left \vert \sum \limits_{j=m+1}^{\infty} \left(\frac{\lambda_j}{\lambda_{*}}\right)^{n}\alpha_j e_j \right \vert \right \vert^2_{L^2(\Omega, c^{-2}dx)} &= \sum \limits_{j=m+1}^{\infty}\left \vert \left \vert \left(\frac{\lambda_j}{\lambda_{*}}\right)^{n}\alpha_j e_j \right \vert \right \vert^2_{L^2(\Omega, c^{-2}dx)}\\
&\leq \left(\frac{\lambda_{m+1}}{\lambda_{*}}\right)^{n} \sum \limits_{j=m+1}^{\infty}\alpha_j \rightarrow 0,
\end{eqnarray*}
as $n \rightarrow \infty$.
Thus
\[\int_{\Omega} g c_1^{-2}\varphi dx=0, \ \ \hbox{for every harmonic functions} \ \ \varphi,\]
where $g$ satisfies $L(g)=\Delta^{-1}(c_1^{-2}g)=\lambda_*g$. By Lemma \ref{LemmaNoClusterCriticalPoints} we have $g= \sum \limits_{j=1}^m \alpha_j e_j \equiv 0$ in $\Omega $, which is a contradiction. Thus $F_0\equiv 0$ and the proof is complete.  \hfill $\Box$ \\ \\ 
Note that Proposition \ref{SpectrumProp} also implies  $\mathcal{H}=\{0\}$, where $\mathcal{H}$ is the Hilbert space defined in the statement of Lemma \ref{Hilbertspace}. \\ \\
{\bf Proof of Theorem \ref{sourceUniqueness}.} The proof follows directly from Theorem \ref{PatternTheo} and Proposition \ref{SpectrumProp}. \hfill $\Box$\\ \\
$\mathbf{Acknowledgement.}$  The authors are grateful to Professor Plamen Stefanov and Gunther Uhlmann for several helpful comments that improved the results of this paper.  Amir Moradifam is supported by NSF grant DMS-1715850. 

\nocite{*}
\bibliography{PATref}{}
\bibliographystyle{plain}

 \end{document}